\documentclass{conm-p-l}

\copyrightinfo{2007}{}

\usepackage{graphicx}

\newtheorem{theorem}{Theorem}[section]
\newtheorem{lemma}[theorem]{Lemma}

\theoremstyle{definition}

\theoremstyle{remark}
\newtheorem{remark}[theorem]{Remark}

\numberwithin{equation}{section}

\newcommand{\tu}{\tilde{u}}

\newcommand{\e}{\epsilon}

\newcommand{\dl}{\delta}

\newcommand{\pa}{\partial}

\newcommand{\om}{\omega}

\newcommand{\na}{\nabla}

\newcommand{\ZZ}{\mathbb{Z}^2/\{0\}}

\newcommand{\tk}{\tilde{k}}
\newcommand{\hu}{\hat{u}}
\newcommand{\tom}{\tilde{\omega}}

\begin{document}

\title[A Recurrence Theorem]{A Recurrence Theorem on the Solutions to the 2D Euler Equation}

\author{Y. Charles Li}
\address{Department of Mathematics, University of Missouri, 
Columbia, MO 65211, USA}
\curraddr{}
\email{cli@math.missouri.edu}
\thanks{}


\subjclass{Primary 37, 76; Secondary 35, 34}
\date{}

\dedicatory{}


\begin{abstract}
In this article, I will prove a recurrence theorem which says that any $H^s(\mathbb{T}^2)$ ($s>2$)  solution to the 2D Euler equation returns repeatedly to an arbitrarily small $H^0(\mathbb{T}^2)$  neighborhood.
\end{abstract}

\maketitle










\section{Introduction}

In finite dimensions, the Poincar\'e recurrence theorem can be proved from the basic properties 
of a finite measure. In infinite dimensions, it is difficult to establish a natural finite measure, especially by extending a finite dimensional finite measure. A natural alternative is the Banach norm 
which can be viewed as a counterpart of the probability density.  An interesting problem is to 
study the Poincar\'e recurrence problem for 2D Euler equation of fluids \cite{Li07}. 
Nadirashvili \cite{Nad91} gave an example of Poincar\'e non-recurrence near a particular solution 
of the 2D Euler equation defined on an annular domain. The proof in \cite{Nad91} was fixed up in
\cite{Li07}. We believe that the Poincar\'e recurrence will occur more often when the 
2D Euler equation is defined on a periodic domain. The theorem to be proved in this article 
is a result along this line. On a periodic domain, the 2D Euler equation is a more natural Hamiltonian 
system than e.g. on an annular domain. 

One final note is that any solution of the 2D Euler equation defines a non-autonomous integrable 
Hamiltonian two dimensional vector field \cite{LL06}. The trajectories of this vector field are the fluid 
particle trajectories. This is the so-called Lagrangian coordinates. The integrability was proved 
in the usual extended coordinates of two spatial coordinates, the stream function, and an extra 
temporal variable. Due to the extra temporal variable, most of the invariant subsets are of infinite
volume. Only in special cases e.g. the solution of the 2D Euler equation is periodic or 
quasi-periodic in time,  one can find invariant subsets of finite volume. This indicates that in the 
Lagrangian coordinates, recurrence is a rare event \cite{Shn97}.

\section{Finite Dimensional System}

In finite dimensions, our theorem is a result of a simple compactness argument. But it states some 
interesting fact. In infinite dimensions, the compactness argument is more complicated.
\begin{theorem}
Let $f \ : \  \mathbb{R}^n \mapsto \mathbb{R}^n$ be a map and $A$ be a compact invariant subset.
Then for any $x \in A$ and any $\dl >0$, there is a $x_* \in A$ such that
\[
f^{m_j}(x) \in B_\dl (x_*)
\]
where $\{ m_j \}$ is an infinite sequence of positive integers and $B_\dl (x_*)$ is the open ball of 
radius $\dl$ centered at $x_*$.
\end{theorem}
\begin{proof}
It is clear that $\{ B_\dl (y)\}_{y\in A}$ is an open cover of $A$, thus there is a finite subcover 
$\{ B_\dl (y_k)\}_{k=1,\cdots , K}$. For any $m=0,1,2,\cdots$; $f^m(x) \in A$; thus $f^m(x) \in 
B_\dl (y_k)$ for some $k$. Therefore, there is at least one $k$ such that an infinite subsequence of 
$\{ f^m(x) \}$ is included in $B_\dl (y_k)$. The theorem is proved.
\end{proof}
\begin{remark}
Of course the theorem is still true when replacing $\mathbb{R}^n$ by a Banach space or a topological
space. But compactness is a very restricted concept in infinite dimensions.
\end{remark}

\section{2D Euler Equation}

The 2D Euler equation
\[
\pa_t u+(u\cdot \na )u = -\na p, \quad \na \cdot u = 0
\]
is globally well-posed in $H^s(\mathbb{T}^2)$ ($s>2$) where $\mathbb{T}^2$ is the 2-tori. We also
require that 
\[
\int_{\mathbb{T}^2} u \ dx = 0.
\]
Denote by $\om$ the vorticity, $\om = \pa_1u_2-\pa_2u_1$. A well-known fact is the following lemma.
\begin{lemma}
\begin{equation}
\int_{\mathbb{T}^2} |\na u|^2 dx = \int_{\mathbb{T}^2} \om^2 dx. 
\label{vde}
\end{equation}
\end{lemma}
\begin{proof}
By the incompressibility condition,
\[
\int_{\mathbb{T}^2}\left [ (\pa_1u_1)^2+(\pa_2u_2)^2\right ] dx = -2 \int_{\mathbb{T}^2}
(\pa_1u_1)(\pa_2u_2) dx.
\]
Since
\begin{eqnarray*}
& & \int_{\mathbb{T}^2}\left [ (\pa_1u_1)(\pa_2u_2) - (\pa_2u_1)(\pa_1u_2)\right ] dx  \\
& & = \int_{\mathbb{T}^2}\left [\pa_1(u_1\pa_2u_2) - \pa_2(u_1\pa_1u_2)\right ] dx = 0 ,
\end{eqnarray*}
we have
\begin{eqnarray*}
\int_{\mathbb{T}^2} |\na u|^2 dx &=& \int_{\mathbb{T}^2}\left [ (\pa_1u_1)^2+(\pa_2u_1)^2
+ (\pa_1u_2)^2+(\pa_2u_2)^2\right ] dx  \\
&=& \int_{\mathbb{T}^2}\left [ (\pa_2u_1)^2+(\pa_1u_2)^2-2(\pa_1u_1)(\pa_2u_2)\right ] dx  \\
&=& \int_{\mathbb{T}^2}\left [ (\pa_2u_1)^2+(\pa_1u_2)^2-2(\pa_2u_1)(\pa_1u_2)\right ] dx  \\
&=& \int_{\mathbb{T}^2}\om^2 dx. 
\end{eqnarray*}
\end{proof}
\begin{remark}
Using Fourier series, one can prove this lemma by direct calculation:
\begin{eqnarray*}
u(k) &=& (ik_2, -ik_1) \frac{1}{|k|^2}\om (k), \\
\int_{\mathbb{T}^2} |\na u|^2 dx &=& \sum_{k \in \ZZ}[k_1^2k_2^2+k_2^4+k_1^4+k_1^2k_2^2]
\frac{1}{|k|^4}\om (k)^2 \\
&=& \sum_{k \in \ZZ}\om (k)^2 = \int_{\mathbb{T}^2} \om^2 dx. \quad \Box
\end{eqnarray*}
\end{remark}
Notice that the kinetic energy
\[
E= \int_{\mathbb{T}^2} |u|^2 dx
\]
and the enstrophy
\[
G= \int_{\mathbb{T}^2} \om^2 dx
\]
are two invariants of the 2D Euler flow.
\begin{lemma}
For any $C>0$, the set 
\[
S=\left \{ u \ \bigg | \  \int_{\mathbb{T}^2} \om^2 dx \leq C \right \}
\]
is compactly embedded in $L^2(\mathbb{T}^2)$ of $u$. That is, the closure of $S$ in 
$L^2(\mathbb{T}^2)$ is a compact subset of $L^2(\mathbb{T}^2)$ .
\label{compe}
\end{lemma}
\begin{remark}
In a simpler language, an enstrophy ball is compactly embedded in the kinetic energy space.
This lemma is the well-known Rellich lemma. In the $\mathbb{T}^2$ setting, we will present 
the proof. There are many versions of the lemma and its proof. We follow that of \cite{Fol76}.
$\quad \Box$
\end{remark}
\begin{proof}
By Lemma \ref{vde},
\[
 \int_{\mathbb{T}^2} |\na u|^2 dx \leq C.
 \]
 Let $\{ u^{(j)} \}$ be a sequence in $S$. For any $\e >0$, choose $K>0$ such that 
\[
2K^{-2}C< \e /2.
\]
For any $k \in \ZZ$, $|k| \leq K$, the Fourier coefficients $\{ u^{(j)}(k) \}$ is a bounded set (e.g. 
bounded by $\sqrt{C}$). Therefore there is a convergent subsequence $\{ u^{(m_j)}(k) \}$. For 
a different $\tk$, $\{ u^{(m_j)}(\tk) \}$ is a bounded set again, thus there is a further convergent subsequence. Iterating on all such $k$ ($|k| \leq K$ finitely many), one can find a subsequence
$\{ u^{(n_j)} \}$ such that $\{ u^{(n_j)}(k) \}$ is uniformly convergent for $|k| \leq K$. This is a simpler 
version of the usual diagonal argument. Next we show that $\{ u^{(n_j)} \}$ is a Cauchy sequence in 
$L^2(\mathbb{T}^2)$ .
\begin{eqnarray*}
\| u^{(n_j)} - u^{(n_\ell )} \|_{L^2(\mathbb{T}^2)} &=& \sum_{|k| \leq K, k\neq 0} |u^{(n_j)}(k) - 
u^{(n_\ell )}(k)|^2 \\
& & + \sum_{|k| > K} |u^{(n_j)}(k) - u^{(n_\ell )}(k)|^2 \\
&\leq&  \sum_{|k| \leq K, k\neq 0} |u^{(n_j)}(k) - u^{(n_\ell )}(k)|^2 \\
& & + K^{-2}\sum_{|k| > K} |k|^2|u^{(n_j)}(k) - u^{(n_\ell )}(k)|^2.
\end{eqnarray*}
The second term is less than $2K^{-2}C< \e /2$. The first term is less than $\e /2$ when $j$ and 
$\ell$ are sufficiently large since $\{ u^{(n_j)}(k) \}$ is uniformly convergent for $|k| \leq K$. So 
$\{ u^{(n_j)} \}$ is a Cauchy sequence in $L^2(\mathbb{T}^2)$, and is convergent inside the 
closure of $S$ in $L^2(\mathbb{T}^2)$. Let $\{ v^{(j)} \}$ be a sequence of the accumulation 
points of $S$ in $L^2(\mathbb{T}^2)$. Then we can find a sequence $\{ u^{(j)} \}$ in $S$ such that 
\[
\| v^{(j)}-u^{(j)} \|_{L^2(\mathbb{T}^2)} < 1/j.
\]
Let $\{ u^{(n_j)} \}$ be the convergent subsequence, then $\{ v^{(n_j)} \}$ is also a convergent subsequence. Thus the closure of $S$ in $L^2(\mathbb{T}^2)$ is a compact subset of 
$L^2(\mathbb{T}^2)$.
\end{proof}

\begin{theorem}
For any $\tu \in H^s(\mathbb{T}^2)$ ($s>2$), any $\dl >0$, and any $T>0$; there is a $u^* 
\in H^s(\mathbb{T}^2)$ such that
\[
F^{m_jT}(\tu ) \in B^0_\dl (u^*)=\{ \hu \in H^s(\mathbb{T}^2) \  | \  \| \hu -u^* \|_{H^0(\mathbb{T}^2)}
< \dl \}
\]
where $\{ m_j \}$ is an infinite sequence of positive integers, and $F^t$ is the evolution operator 
of the 2D Euler equation.
\end{theorem}
\begin{proof}
Choose the $C$ in Lemma \ref{compe} to be
\[
2\int_{\mathbb{T}^2} |\na \tu |^2 dx = 2\int_{\mathbb{T}^2} \tom^2 dx.
\]
Define two sets:
\begin{eqnarray*}
S &=& \left \{ u \  \bigg | \  \int_{\mathbb{T}^2} \om^2 dx \leq 2 \int_{\mathbb{T}^2} \tom^2 dx 
\right \} , \\
S_1 &=& \left \{ u \in H^s(\mathbb{T}^2) \  \bigg | \  \int_{\mathbb{T}^2} \om^2 dx \leq 2
 \int_{\mathbb{T}^2} \tom^2 dx \right \} .
 \end{eqnarray*}
Notice that $S_1$ is invariant under the 2D Euler flow, and $S_1$ is a dense subset of $S$,
$S_1 = S \cap H^s(\mathbb{T}^2)$. By Lemma \ref{compe}, the closure of $S$ in $L^2(\mathbb{T}^2)
=H^0(\mathbb{T}^2)$ is a compact subset. For any $u \in S$, denote by
\[
B_{\dl /2}(u)=\{ u \in H^0(\mathbb{T}^2) \  | \  \| v -u \|_{H^0(\mathbb{T}^2)}
< \dl /2\} .
\]
All these balls $\{ B_{\dl /2}(u) \}_{u \in S}$ is an open cover of the closure of $S$ in 
$H^0(\mathbb{T}^2)$. Thus there are finitely many $u_1, \cdots ,u_N \in S$ and 
$\{ B_{\dl /2}(u_n) \}_{n=1, \cdots , N}$ is a finite cover. Since $S_1$ is dense in $S$, for each such 
$u_n$, one can find a $u^*_n \in S_1$ such that 
\[
\| u_n - u^*_n \|_{H^0(\mathbb{T}^2)} \leq \| u_n - u^*_n \|_{W^{1,2}(\mathbb{T}^2)}  < \dl /4,
\]
where $W^{1,2}(\mathbb{T}^2)$ is the Sobolev space with norm $\| u\|^2 = \int_{\mathbb{T}^2}
|\na u|^2 dx$. All the balls
\[
B_{\dl }(u^*_n)=\{ v \in H^0(\mathbb{T}^2) \  | \  \| v -u^*_n \|_{H^0(\mathbb{T}^2)} < \dl \}
\]
still covers $S$, thus covers $S_1 = S \cap H^s(\mathbb{T}^2)$. Let $B^0_{\dl }(u^*_n)=
B_{\dl }(u^*_n) \cap H^s(\mathbb{T}^2)$,
\[
B^0_{\dl }(u^*_n)= \{ \hu \in H^s(\mathbb{T}^2) \ | \   \| \hu - u^*_n \|_{H^0(\mathbb{T}^2)} < \dl \} .
\]
Then 
\[
S_1 \subset \bigcup_{n=1}^N B^0_{\dl }(u^*_n) .
\]
By the invariance of $S_1$ under the 2D Euler flow $F^t$. There is at least one $n$ such that 
an infinite subsequence of $\{ F^{mT}(\tu )\}_{m=0,1, \cdots }$ is included in $B^0_{\dl }(u^*_n) $.
The theorem is proved.
\end{proof}

\end{document}